\documentclass[10pt,final]{amsart}
\usepackage[latin1]{inputenc}
\usepackage{latexsym}
\usepackage{a4wide}
\usepackage{amscd}
\usepackage{graphics}
\usepackage{amsmath}
\usepackage{amssymb}
\usepackage{mathrsfs}
\input xy
\xyoption{all}
\usepackage{bm}

\newcommand{\N}{\mathbb{N}}                    % the natural numbers
\newcommand{\Z}{\mathbb{Z}}                    % the integer numbers
\newcommand{\R}{\mathbb{R}}                    % the real line
\newcommand{\Ddtt}{\tfrac{ D^2}{ dt^2}}
\newcommand{\finedim}{\hfill $\Box$\\}          % end of proof
\newcommand{\sgn}{\mathrm{sgn\,}}              % signum
\newcommand{\diag}{\mathrm{diag\,}}            % diagonal matrix
\newcommand{\poset}{{\preceq^{\ominus}}}

\newcommand{\Id}{\mathrm{Id}}

\newcommand{\Spl}{\mathrm{Sp}}

\newcommand{\n}[1]{{\bf #1}}
\newcommand{\Bsym}{\mathrm B_{\mathrm{sym}}}

\newcommand{\pparagraph}[1]{{\bf{#1}}}
\def\virgA{``}

\newtheorem{mainthm}{\sc Theorem}          % numbered absolutely
\newtheorem{thm}{\sc Theorem}[section]      % numbered within each section
\newtheorem{cor}[thm]{\sc Corollary}        % numbered along with Theorem
\newtheorem{lem}[thm]{\sc Lemma}            % numbered along with Theorem
\newtheorem{prop}[thm]{\sc Proposition}    % numbered along with Theorem
\newtheorem{defin}[thm]{\sc Definition}    % numbered along with Theorem
\newtheorem{rem}[thm]{\sc Remark}          % numbered along with Theorem
\title[Asymptotically linear Hamiltonian systems]{A Multiplicity result
for a class of strongly
indefinite asymptotically linear second order systems}
\author[A. Capietto, F. Dalbono and A. Portaluri]{Anna Capietto, Francesca
Dalbono and Alessandro Portaluri}
\address{{\sl Anna Capietto {\rm and} Francesca Dalbono:} \hfill\break\indent Dipartimento di Matematica,\hfill\break\indent
Universit\`a di Torino, \hfill\break\indent Via Carlo Alberto 10,
Torino \hfill\break\indent Italy}   \email{anna.capietto@unito.it}
\email{francesca.dalbono@unito.it  \hfill\break\indent tel.(39)0116702914} 

\address{ {\sl Alessandro Portaluri:}  \hfill\break\indent
Dipartimento di Matematica \virgA Ennio De Giorgi",\hfill\break\indent
Ex-collegio Fiorini,  Universit\`a del Salento, \hfill\break\indent
Via per Arnesano, Lecce \hfill\break\indent Italy }
\email{alessandro.portaluri@unile.it  \hfill\break\indent tel.(39)0832297589}
 
\thanks{%
\noindent The work of the first two authors has been performed in the frame of the Italian M.I.U.R. project "Variational and topological methods in the study of nonlinear phenomena" and of the Italian GNAMPA-I.N.d.A.M. project "Nonlinear analysis techniques for boundary value problems associated to differential equations". The third author was partially supported by Italian M.I.U.R. project
"Variational Methods and Nonlinear Differential Equations".}

\begin{document}

\begin{abstract}
We prove a multiplicity result for a class
of strongly indefinite nonlinear second order asymptotically linear systems
with Dirichlet boundary conditions. The key idea for the proof is to
bring together the classical shooting method and the Maslov index of the linear Hamiltonian 
systems associated to the asymptotic limits of the given nonlinearity.
\end{abstract}
\maketitle

\section{Introduction}

\noindent In this paper we deal with a second order nonlinear boundary value problem of the form
\begin{equation}
\label{eq:Iprin}
\left\{\begin{array}{ll}
J u''(t) + S(t,u(t))u(t) =0\\
u(0)=0=u(1),
\end{array}\right.
\end{equation}
where 
\begin{equation}\label{eq:simmetriaintro}
J=\begin{pmatrix} \Id_{n-\nu} & 0 \\
0 & -\Id_\nu
\end{pmatrix}
\end{equation}
and $S:[0,1]\times \R^n \to \Bsym(\R^n)$ is continuous. It is useful to observe that any
gradient system of the form
\begin{equation}\label{eq:IsistemanonlinearconVintro}
J u''(t) + \nabla V(t,u(t)) =0
\end{equation}
with $\nabla V(t,0)=0$ is of the form $J u''(t) + S(t,u(t))u(t) =0$. Indeed, it is sufficient to set  
\begin{equation*}
S(t, u) := \int_0^1 D^2 V(t, su) ds, \qquad \forall\, (t,u) \in [0,1]
\times \R^n.
\end{equation*} 

\noindent We are concerned with the existence and multiplicity of solutions to \eqref{eq:Iprin} with prescribed nodal properties.

\noindent Before describing the set of assumptions on the nonlinearity and the method of the proofs and developing some remarks with the literature, we focus on some motivations arising from the study of differential equations on manifolds. More precisely, we wish to (very briefly) describe the question of the study of the conjugate points along a perturbed geodesic in a semi-Riemannian manifold in relation with systems of the form $J u''(t) + S(t)u(t) =0$. For a general reference, we quote the book \cite{BeeEhrEas}; recent results can be found, among others, in 
\cite{Hel94},\cite{MPP},\cite{MusPejPor07},\cite{PicPorTau04},\cite{PorEJDE}. See also \cite{fpr}.

\noindent Let $M$ be a
smooth semi-Riemannian manifold, i.e. a $C^\infty$, $n$-dimensional
manifold $M$ endowed with a (semi)-Riemannian  metric (i.e. a
non-degenerate symmetric two-form $g$ of constant index $\nu\in \{0,
\dots, n\}$). Denoted by $D$ and by  $\frac{D}{dt}$ respectively the
associated Levi-Civita connection and the covariant derivative of a
vector field along a smooth curve $\gamma$, a {\em
perturbed geodesic\/} or briefly a {\em p-geodesic\/} is a smooth
curve $\gamma:[0,1] \rightarrow M$ which satisfies the differential
equation
\begin{equation}\label{2.2}
\dfrac{D}{dt} \gamma '(t) +\nabla_g V(t,\gamma(t))=0
\end{equation}
where  $\nabla_g V$ denotes gradient of $V (t,-) $ with respect to
the metric tensor $g.$

\noindent Let $\gamma$ be a $p$-geodesic between two fixed points $p,q \in M$. A
vector field $\xi$ along $\gamma$  is called a {\em Jacobi field\/}
if it verifies the  linear second order differential equation
\begin{equation}\label{2.5}
\Ddtt\xi (t) + R(\gamma ' (t) , \xi (t)) \gamma '(t) + D_{\xi
(t)}\nabla V(t, \gamma(t))=0,
\end{equation}
where $R$ is the curvature tensor of $D.$ Given a $p$-geodesic
$\gamma$, an instant $t \in (0,1]$ is said to be a {\em conjugate
instant\/} if there exists at least one non zero Jacobi field with
$\xi(0)=\xi (t)=0.$ The corresponding  point $ q=\gamma(t)$ on $M$
is said to be a {\em conjugate point\/} to the point $p =\gamma(0)$
along $\gamma.$ 

\noindent Equation \eqref{2.5} can be then written in the form $J u''(t) + S(t)u(t) =0$. Indeed (we refer for more details to \cite{MPP}), given a perturbed geodesic $\gamma$, a vector field $\xi$ along $\gamma$ can be written as $\xi(t)=\sum_{i=1}^{n}u_i(t)e^i(t),$ being 
$\{e^1, \dots, e^n\}$ a $g$-frame along $\gamma$ (this means that the $e^i$s are pointwise $g$-orthogonal and $g(e^i(t),e^i(t))=\epsilon_i,$ with $\epsilon_i=1$ for $i=1,\dots, n-\nu$ and $\epsilon_i=-1$ for $i=n-\nu+1,\dots, n,$ for all $t$). Equation \eqref{2.5} can be thus transformed into the second order system of the form:
\begin{equation}
\epsilon_i u''_i(t)+\sum_{j=1}^{n}S_{ij}(t)u_j(t)=0, \qquad i=1,\dots,n,
\end{equation}
where $S_{ij}=g(R(\gamma',e^i)\gamma'+D_{e^i}\nabla V(\cdot,\gamma),e^j)$.

\noindent Here we deal with a class of nonlinearities $S$ which have an "asymptotically linear" behaviour at zero and infinity; more precisely, we assume
\begin{enumerate}
\item[$(V_0)$] $S(t, \xi) \xi  = S_0(t) \xi + o(|\xi|)$ for $|\xi|\to 0$, uniformly in
$t$;
\item[$(V_\infty)$] $S(t, \xi)\xi= S_\infty(t) \xi + o(|\xi|)$ for $|\xi|\to +
\infty$, uniformly in $t$.
\end{enumerate}
In the study of nonlinear boundary value problems, it is quite frequent to meet this set of hypotheses; they give rise, in the suitable context, to two indices which contain some information on the behaviour of the problem at zero and infinity. Then, roughly speaking, the bigger the gap between these two indices the greater the number of solutions of the nonlinear BVP. As for BVPs with separated boundary conditions, we refer, among others, to  \cite{CapDal}, \cite{CapDam}, \cite{CapDamPap}, \cite{DalReb08}, \cite{DalZan07}, \cite{Do05}, \cite{Fo-96}, \cite{liu}, \cite{P}, \cite{Snla}, \cite{Z}. For results in the same spirit but for the periodic problem we refer to the pioneering work of Amann-Zehnder \cite{AZ}  and to the more recent works \cite{BeFo-94}, \cite{F},  \cite{Izy99}, \cite{LSW}, \cite{mareza}, \cite{szul}. We take as a starting point of the present research the paper \cite{CapDamPap}, where it is treated the particular case $\nu=0$. By classical shooting methods, when $n=1$ (a scalar equation) and  $\nu=0$ (the positive definite case) a gap condition expressed in terms of the rotation number is sufficient for the existence of multiple solutions with a prescribed number of zeros. However, also in the positive definite case, when $n \geq 2$ more assumptions on $S$ (together with the use of the Maslov index) are needed (cf. Remark 4.10 in \cite{CapDamPap}). More problems arise in case $\nu \neq 0$. In this paper, in order to deal with {\it systems} we require (on the lines of \cite{CapDamPap}) a diagonality condition (cf. $(V_3)$) on the restriction of $S$ to the $(n-1)$ coordinate hyperplanes of $\R^n$; the difficulty due to the {\sl  indefiniteness} of $J$ is treated by assuming that $S$ is $J$-commuting (the "split" condition $(V_1)$). In order to distinguish the solutions we use a generalized shooting method initiated in the linear case by L. Greenberg in \cite{greenberg} (similar ideas can be found also in \cite{Re-80}). Indeed, 
taking advantage of the symplectic  
structure of the linear Hamiltonian system associated to a linear second order system of the form \eqref{eq:Iprin}, L. Greenberg \cite{greenberg} has generalized the  
well-known concepts of elementary phase-plane analysis; he has developed the concepts of lagrangian plane, phase angle and crossing, which correspond, in the case of planar systems, to the  
notions of line through the origin, polar coordinate and  
zero of a real function, respectively. Thus, in the present paper we adapt the notion of $h$-type solution given in \cite{CapDamPap} (inspired by \cite{greenberg}) to our more general (indefinite) framework. 

\noindent For the statement of our main results (Theorems 1, 2, 3), we focus on the Maslov indices $m_0$ and $m_{\infty}$ of the linear systems $Ju''+S_0(t)u=0$ and $Ju''+S_{\infty}(t)u=0,$ respectively. By these indices, we define (cf. \eqref{eq:insiemeTNOnpositividefinito}) a set $\mathscr T$ whose non-emptiness is a sufficient condition (Theorem 1) for the existence of multiple solutions to the given BVP. Theorems 2 and 3 provide sufficient conditions for the non-emptiness of $\mathscr T$ (cf. also Proposition \ref{asimmetria} and Proposition \ref{sceltaalphainfty}). For the proofs, we employ the concept of phase angle and the generalized shooting method in order to reformulate the problem in terms of the existence of zeros  
of an $N$-dimensional vector field. Then, we apply a version of the Miranda fixed point theorem given in \cite{CapDamPap}. 

\noindent  Our results extend the main result in \cite{CapDamPap} to the case $\nu \neq 0$. In particular, they represent a generalization to indefinite systems of one-dimensional results (we refer, among others, to  \cite{DalZan07}, \cite{DinSan}, \cite{Est}, \cite{GHU}, \cite{Snla}) where multiplicity is obtained through a comparison between the behaviour of the nonlinearity at zero and infinity. It is worth noticing (as it is explained in detail in Remark \ref{notcomp}) that our approach is similar to the one in \cite{CapDamPap}, but our results are not a consequence of Theorem 4.7 in \cite{CapDamPap}.

\noindent Some final comments on the literature are in order. Indeed, we wish to point out that multiplicity for asymptotically linear problems has been achieved, with various different methods, by some authors without a  diagonality condition of the form $(V_3)$; however, other restrictions are needed. For a comprenhensive reference, we refer to the book by J.Mawhin-M.Willem \cite{MaWi}. In particular, in \cite{Do05}, \cite{P} and \cite{Z} (the second and third in the framework of elliptic PDEs) the potential $V$ in \eqref{eq:IsistemanonlinearconVintro} is even; in our work, we do not deal with such kind of restriction (cf. Proposition \ref{asimmetria} for more details). It is also worth noticing that no diagonality condition is imposed also  in \cite{DalReb08}, where, in turn, it is required (for the planar case) a sign condition on $S$ and solutions which are not of $h$-type are obtained. Last but not least, it is important to remark that by using the Maslov index we have been able to treat indefinite problems, which may have an infinite Morse index (and whose related action functional is unbounded).

\noindent We end this Section with a list of notations. 

\noindent In what follows, we set $\N^{\ast} := \N \setminus \{0\}$ and ${\bf n}:=\{i \in \N: i\in [1,n]\}$. We denote by $\underline 1$
the vector whose components are equal to $1$ and 
we write $\underline j$ for the vector whose components are the elements on the
diagonal of the matrix $J$. The $(n \times n)$ identity matrix is denoted by $\Id_n$. Given the posets $(Z_1,\preceq)$ and $(Z_2, \preceq^{op})$,
where $\preceq^{op}$ denotes its dual order, we define the poset
$(Z_1 \oplus Z_2,\preceq^\ominus)$ as the direct sum of the above
defined posets; we also define the poset $(Z_1 \oplus Z_2,\preceq^\oplus)$ as the direct sum of
the posets $(Z_1,\preceq)$ and $(Z_2, \preceq)$. By $\langle \cdot, \cdot \rangle$ we mean the scalar product. We denote by $\Bsym(\R^n)$ and $\Spl(n)$ the set of $(n \times n)$ symmetric and symplectic matrices, respectively. Consider the $(2n\times 2n)$ matrix
\begin{equation} \label{eq:sigma} 
\sigma_n =\begin{pmatrix} 0 & \Id_n
\cr -\Id_n & 0
\end{pmatrix},
\end{equation}
and the standard symplectic form $\omega(z_1, z_2) := \langle \sigma_n z_1, z_2 \rangle$ for $z_i \in \R^{2n}$. We denote by
$\mathscr L:=\mathscr L(\R^{2n}, \omega)$ the set of all Lagrangian subspaces
of the symplectic space $(\R^{2n}, \omega)$. For $l_0\in \mathscr L$, the set
$\Sigma(l_0) = \{ l\in \mathscr L: l \cap l_0 \not= \{0\}\}$ is the
{\em train\/} or the {\em Maslov cycle\/} of $l_0$.  The Lagrangian subspace $\{0\}\oplus {\R}^n \subset {\R}^{2n}$ will be denoted by $L_0$ and we refer to it with the name of {\it vertical Lagrangian.}

\bigskip
\bigskip

\section{Linear symplectic preliminaries}
\label{lin} 
\noindent In this section we first recall (according to \cite{RobSal93}) the definition of Maslov index; then, we introduce the notion of phase angle (on the lines of \cite{greenberg}). Finally, we state and prove a Sturm-type theorem.

\vspace{3mm}
 
\noindent \pparagraph{The Maslov index.} The
Maslov index is a semi-integer homotopy invariant with fixed endpoints of paths $l$ of
Lagrangian subspaces of the symplectic vector space $(\R^{2n}, \omega)$
which gives an algebraic count of non transverse intersections of
the family $\{l(t)\}_{t \in [0,1]}$ with a given Lagrangian subspace
$l_{0}$. For each $C^1$-curve $l \colon [0,1] \to \mathscr L$ of Lagrangian
subspaces, we say that $t_0 \in [0,1]$ is a 
{\em crossing\/} for $l$ if $l(t_0)\cap l_0\not=\{0\}$,
i.e. $l(t_0)\in\Sigma(l_0)$. Consider $t_0 \in [0,1]$ and denote by $W$ a lagrangian complement of $l(t_0)$. For $z \in l(t_0)$ and $t$ in a neighbourhood of $t_0$,  let $w(t)$ be the unique vector s.t.
\[
w(t) \in W, \qquad z+w(t) \in l(t).
\] 
Let us then set, for $z \in l(t_0)$, 
\[
Q(l,t_0)(z)= \frac{d}{dt} \omega(z,w(t)) |_{t=t_0}.
\]
Note that $Q$ is independent of $W$. At each crossing $t_0 \in [0,1]$ we define the {\em
crossing form\/} $\Gamma$ as the quadratic form
\[
\Gamma (l, l_0, t_0)\, :=\, Q(l, t_0) |_{l(t_0) \cap l_0}.
\]

\noindent The crossing $t$ is {\em regular\/} if the
crossing form is nonsingular.  It is easy to check that regular crossings are isolated and
therefore on a compact interval they are in a finite number. Assuming that $l$ has only regular crossings, we can give the following
\begin{defin}\label{def:RobSalindex}
The Maslov index of the Lagrangian path $l$ relative to the
Lagrangian subspace $l_0$ is the semi-integer defined by
\begin{equation}\label{eq:MaslovdiRobSal}
\mu_{l_0}(l, [0,1])\, := \, \dfrac12\,\sgn\, \Gamma (l, l_0, 0)\, +
\,\sum_{t\in\left]0,1\right[}\sgn\, \Gamma (l, l_0, t)\, +
\,\dfrac12\,\sgn\, \Gamma (l, l_0, 1),
\end{equation}
where $\sgn$ denotes the signature of a quadratic form and the summation runs over all crossings $t$.
\end{defin}

\noindent  It is a standard fact that the above definition can be extended to the case when there exist non-regular crossings. The above definition comes from \cite{RobSal93}, to which we refer for all the details; we just note that the construction in \cite{RobSal93} is developed on the basis of the form $-\sigma_n$ instead of $\sigma_n$. For other useful references on the Maslov index, we refer to \cite{abbo}, \cite{A},\cite{long}. 

\noindent Given $\psi \colon [0,1] \to \Spl(2n)$ a continuous path of
symplectic matrices and $l_0$ a lagrangian subspace, then we define the
{\em Maslov index\/} of $\psi$ as follows
\[
\mu_{l_0}(\psi, [0,1])\, := \, \mu_{l_0}(\psi (\cdot) (l_0),[0,1]).
\]

\noindent We finally give the following

\begin{defin}\label{def:multiplicitycrossing}
For $k \in \n{n}$,  we say that a crossing $t_0 \in [0,1]$ has {\em multiplicity\/}  $k$ if $k$ if the dimension of the intersection between $l_0$ and
$l(t_0)$.
\end{defin}

\vspace{3mm}

\noindent \pparagraph{A phase angle analysis.}\label{sec:phaseangles} Let
$S\colon [0,1] \to \Bsym(\R^n)$ be a continuous path of symmetric matrices and let us
consider the linear second order system
\begin{equation}\label{eq:linearesecondoordine}
Ju''(t) + S(t)u(t) =0,  \qquad t \, \in [0,1],
\end{equation}
where $J$ is defined in \eqref{eq:simmetriaintro}. By performing the
change of coordinates $v= Ju'$, \eqref{eq:linearesecondoordine} can be written as the following first order system
\begin{equation*}
\left\{\begin{array}{ll} u'(t) = Jv(t) \\ 
v'(t) = - S (t) u(t).
\end{array}\right.
\end{equation*}
By taking $w=(u,v)$, \eqref{eq:linearesecondoordine} takes the Hamiltonian form
\begin{equation}\label{eq:ham}
w'(t) = \sigma_n H(t) w(t),
\end{equation}
where
\begin{equation}\label{eq:accat}
H(t) = \begin{pmatrix}  S(t) & 0 \cr 0 &  J
\end{pmatrix}
\end{equation}
for each $t \in [0,1]$.
Under this change of variables the Dirichlet boundary conditions become
\[
(w(0), w(1)) \in L_0 \times L_0.
\]

\noindent In what follows, we shall need the following 

\begin{defin}\label{def:splitmatrix}
A {\em split matrix\/} is any matrix commuting with $J$.
\end{defin}
\noindent It is easy to see that any split matrix has the form:
\begin{equation}
S =\begin{pmatrix} A & 0 \\
0 & B
\end{pmatrix},
\end{equation}
where $A$ and $B$ are symmetric $(n-\nu \times n-\nu)$ and $(\nu
\times \nu)$ matrices, respectively. From now on, we will assume the
following condition:
\begin{enumerate}
{\em \item[(S)] the continuous path of symmetric matrices $S\colon
[0,1]\to \Bsym(\R^n) $ is split.\/}
\end{enumerate}
Let $K$ be the continuous path defined pointwise by
\[
K(t):=\sigma_n H(t)= \begin{pmatrix} 0 & 0 & \Id_{n-\nu}& 0 \\
0 & 0 & 0& -\Id_\nu\\
-A_{n-\nu}(t) & 0 & 0 & 0\\
0 & -B_\nu(t) & 0 & 0
\end{pmatrix} \quad \textrm{for} \ \ S(t):= \begin{pmatrix} A_{n-\nu}(t) & 0\\
0 & B_{\nu}(t)
\end{pmatrix}
\]
and let $U$ be the block matrix
\begin{equation}\label{eq:cambiocartalagrangiani}
U:=\begin{pmatrix} \Id_{n-\nu} & 0 & 0
& 0\\
0& 0 & \Id_{n-\nu} & 0\\
0& \Id_{\nu} & 0 & 0\\
0& 0 & 0 &
\Id_{\nu}\\
\end{pmatrix}.
\end{equation}
Denote by $\widetilde K$ the continuous path defined pointwise by
\begin{equation}\label{eq:Kbarra}
\widetilde K(t):=\begin{pmatrix} 0& \Id_{n-\nu}& 0 & 0 \\
-A_{n-\nu}(t) & 0 & 0 & 0 \\
0& 0 & 0 & -\Id_\nu\\
0& 0 & -B_\nu(t) & 0
\end{pmatrix}.
\end{equation}
By the change of coordinates
$\widetilde w=Uw,$ since $ UK=\widetilde K U$, system \eqref{eq:ham} reduces to $\widetilde
w'=\widetilde K \widetilde w$. If $\widetilde w:=(\widetilde w_1,
\widetilde w_2)$, the first order system $\widetilde w'=\widetilde K
\widetilde w$ can be written as
\begin{equation}\label{eq:hamdisaccoppiato}
\left\{\begin{array}{ll}
\widetilde w'_1(t)\ =\  k_A(t) \ \widetilde w_1(t) \\
\widetilde w'_2(t)\ =\ k_B(t)\  \widetilde w_2(t)
\end{array}\right.
\end{equation}
being
\[
k_A(t):= \begin{pmatrix} 0& \Id\\
-A(t) & 0
\end{pmatrix}, \qquad  k_B(t):=\begin{pmatrix} 0& -\Id\\
-B (t)& 0
\end{pmatrix}.
\]
If $L_0^B, L_0^A$ are (respectively) the vertical Lagrangian subspaces
in the symplectic spaces  $(\R^{2\nu}, \sigma_\nu)$ and $(\R^{2(n
-\nu)}, \sigma_{n-\nu})$, we have $UL_0=L_0^A\oplus L_0^B$. Thus by
setting $\widetilde L_0 := L_0^A\oplus L_0^B$, the Dirichlet
boundary conditions can be written as
\[
(\widetilde w(0), \widetilde w(1)) \in \widetilde L_0\times
\widetilde L_0.
\]
For each $t \in [0,1]$, let us consider the $n-\nu$ independent
solutions $\widetilde w_1^1(t), \dots , \widetilde w_1^{n-\nu}(t)$
of
\begin{equation}\label{eq:hamiltnianform1}
\left\{\begin{array}{ll} \widetilde w_1'(t) = k_A(t) \ \widetilde w_1(t)\\%%\qquad t \in[0,1]\\
\widetilde w_1(0)\in L_0^A
\end{array}\right.
\end{equation}
and the $\nu$ independent solutions $\widetilde w^1_2(t), \dots ,
\widetilde w_2^{\nu}(t)$ of
\begin{equation}\label{eq:hamiltnianform2}
\left\{\begin{array}{ll} \widetilde w_2'(t) = k_B(t) \widetilde w_2(t)\\%% \qquad  t \in [0,1]\\
\widetilde w_2(0)\in L_0^B.
\end{array}\right.
\end{equation}
By setting
\[
\widetilde w_1^1(t)=\begin{pmatrix}u_1^1(t)\\
v_1^1 (t)\end{pmatrix} , \dots ,\widetilde
w_1^{n-\nu}(t)=\begin{pmatrix}
u_1^{n-\nu}(t)\\
v_1^{n-\nu}(t)
\end{pmatrix}\]
and
\[
\widetilde w_2^1(t)=\begin{pmatrix}u_2^1(t)\\
v_2^1 (t)\end{pmatrix} , \dots ,\widetilde
w_2^{\nu}(t)=\begin{pmatrix}
u_2^{\nu}(t)\\
v_2^{\nu}(t)
\end{pmatrix},
\]
we can define the two matrices $ X_1(t):=[u_1^1(t), \dots,
u_1^{n-\nu}(t)],$ $ X_2(t):=[u_2^1(t), \dots, u_2^{\nu}(t)]$ and
$X_1'(t):=[v_1^1(t), \dots ,v_1^{n-\nu}(t)]$, $X_2'(t):=[ v_2^1(t),
\dots, v_2^{\nu}(t)]$. Since, for each $t \in [0,1],$ the
vectors $\{\widetilde w_1^1(t), \dots, \widetilde w_1^{n-\nu}(t)\}$
and $\{\widetilde w_2^1(t), \dots ,\widetilde w_2^{\nu}(t)\}$ are
linearly independent, the matrices $\widetilde X_1(t)=[X_1(t)\
X'_1(t)]^T$ and $\widetilde X_2(t)=[X_2(t)\ X'_2(t)]^T$ have rank $n-\nu $ and $\nu,$ respectively; hence the matrix
$X'_j(t)-iX_j(t)$ is invertible for each $j \in \n{2}$. Now we
define, for each $t \in [0,1]$,  the unitary symmetric matrices
\begin{equation}\label{eq:angolifaseblocchi}
Y^j(t):=\big(X'_j(t)+iX_j(t)\big)\big(X'_j(t)-iX_j(t)\big)^{-1}\qquad
\textrm{for}\  j=1,2;
\end{equation}
let us denote by $\{\lambda_l^{j}(t)\}$ their spectrum.
Here we refer, for instance, to \cite[Section 6]{greenberg}. By Kato's
selection Theorem (cf. \cite[Chapter II, Section 6]{Kat80}), for
each $l\in \{1, \dots, n-\nu\}$ there exists a unique continuous map
$\vartheta^1_l : [0,1] \to \R$ such that $\lambda^1_l(t)= e^{2
i\vartheta^1_l(t)}$ with $\vartheta^1_l(0)=0$ and for each $l\in
\{n-\nu+1,\ldots,n\}$ there exists a unique continuous map
$\vartheta^2_l : [0,1] \to \R$ such that $\lambda^2_l(t)= e^{2
i\vartheta^2_l(t)}$ with $\vartheta^2_l(0)=0$.
 
\noindent We are now ready for the following
\begin{defin}
\label{def:angolidifase}
For $j \in \n{2}$, we term {\em phase angles\/} of the system
\eqref{eq:linearesecondoordine} the continuous functions $\theta^j_l\colon [0,1] \to
\R$  obtained by continuously arranging the $n-\nu$ functions
$\vartheta^1_l$ corresponding to $Y^1$ in increasing order and
the $\nu$ functions $\vartheta^2_l$ corresponding to $Y^2$ also
in increasing order. 
\end{defin}

\noindent With this setting, we shall write $\Theta^1=(\theta^1_1, \dots, \theta^1_{n-\nu})$, 
$\Theta^2=(\theta^2_{n-\nu+1}, \dots, \theta^2_n)$ and $\Theta=\Theta^1 \oplus \Theta^2.$
\noindent The following lemma explains the relation between the notions of crossing and of phase angle. To this end, we denote by $\psi$ the fundamental solution associated to the Hamiltonian system
\eqref{eq:hamdisaccoppiato}.

\begin{lem}\label{eq:angolidifasemolteplicita}
The following facts are equivalent:
\begin{enumerate}
\item  $t_0\in [0,1]$ is a crossing for $\psi(\cdot)(\tilde L_0)$ (w.r.t. $\tilde L_0$) of multiplicity
$\mu \in \n{n}$;
\item there exist
exactly $\mu$ different integers $l_1, \dots, l_{\mu} \in \n{n}$ and
there exist $j_1,\ldots,j_\mu\in\n{2}$, $h_1, \dots, h_\mu \in \Z$
such that
\[
\theta^{j_1}_{l_1}(t_0)=  h_1 \pi,\,\, \dots\,\,,
\theta^{j_\mu}_{l_{\mu}}(t_0)= h_\mu \pi.
\]
In particular, $h_k\in\N$ if $j_k=1$ and $h_k\in-\N$ if $j_k=2$.
\end{enumerate}
\end{lem}

\proof This result can be proved by arguing as in \cite[Proposition
3.13]{CapDamPap}. The sign of the phase angles can be easily
deduced by \cite[Lemma 8.2]{greenberg}. \finedim \\

\noindent In what follows we shall write, for each $t \in [0,1],$ the phase angles in
the following form
\begin{equation}\label{eq:angolieinteri}
\theta^1_l(t)= k^1_l(t) \pi + \alpha^1_l(t), \ \  l \in \{1, \dots,
n-\nu\};\qquad \theta^2_l(t)=-k^2_l(t) \pi - \alpha^2_l(t), \ \ l \in
\{n-\nu+1,\ldots,n\},
\end{equation}
where, for each $t \in (0,1]$,  $k^j_l(t) \in \N$ and $
\alpha^j_l(t) \in (0,\pi]$.\\

\noindent Let us now turn to a step which will be crucial for the proof of our main results. \\

\vspace{3mm}

\pparagraph{\noindent Sturm Comparison Principle for scalar equations.} Let us
consider the initial value problem  
\begin{equation}\label{eq:scalar}
\left\{\begin{array}{ll} u''(t)+ a(t)u(t)=0\\
u(0)=0,\ \  u'(0)=1\end{array}\right.
\end{equation} for some
continuous function $a$. In this case, for every $t\in [0,1]$, the
matrix $Y(t)$ reduces to the complex number
\[
Y(t)= \dfrac{\varphi'(t)+ i \varphi(t)}{\varphi'(t)- i \varphi(t)},
\]
where $\varphi$ is the solution of \eqref{eq:scalar}. It is easy to
show that the unique phase angle $\theta(t)$ is the argument, in
polar coordinates, of the complex number $\widetilde\varphi(t) = \varphi'(t) +
i\varphi(t)$. We point out that $\theta(t)$ does not coincide exactly with
the usual polar coordinate in the phase-plane; indeed, the phase
angle $\theta(t)$ is measured in the standard Euclidean plane $(x,y,O)$
starting from the $y$-axis and in the clockwise sense.
\begin{lem}\label{thm:favata}
Consider
\begin{equation}\label{eq:modello}
u''(t)+ a(t)u(t)=0, \qquad t \in [0,1].
\end{equation}
Denoting by $\vartheta^1_a,  \vartheta^2_a$ the phase angles
associated respectively to the Hamiltonian systems:
\[
\left\{ \begin{array}{ll} u'=v_1\\
v_1'=-au
\end{array}\right. \qquad \left\{\begin{array}{ll} u'=-v_2\\
v_2'=au,
\end{array}\right.
\]
we have
\begin{enumerate}
\item $\vartheta^2_a(t)=- \vartheta^1_a(t);$
\item if $a(t) \leq b(t)$ then $\vartheta^1_a(1) \leq \vartheta^1_b(1)$ (or, equivalently, $
\vartheta^2_a(1) \geq \vartheta^2_b(1)$).
 
\end{enumerate}
\end{lem}
\proof We prove (1). By using Pr\"ufer coordinates in the
phase plane $(u, u')$, the following relation holds:
\[
\vartheta^2_a(t)= \arctan\left(\frac{u(t)}{v_2(t)}\right)
=\arctan\left(-\frac{u(t)}{u'(t)}\right)= - \vartheta^1_a(t).
\]
Now the first conclusion in $(2)$ readily follows by
the Sturm comparison principle (cf. \cite{CodLev55}). By (1), we also deduce that if $a(t) \leq b(t)$ then
$\vartheta^2_a(1) \geq \vartheta^2_b(1)$ and this concludes the
proof.\finedim
 \begin{rem}
\label{remrem}
Given $a_i,b_j\in C^0([0,1])$ with $i\in\{1,\ldots,n-\nu\}$,
$j\in\{n-\nu+1,\ldots, n\}$, consider the following uncoupled second order problem
\begin{equation}
\label{eq:diang}
\left\{\begin{array}{ll}
Ju''(t) + \Delta(t) \, u(t) =0,
\qquad t \in [0,1]\\
u(0)=0=u(1),
\end{array}\right.
\end{equation}
where $\Delta(t):=\diag(a_1(t), \dots,a_{n-\nu}(t),b_{n-\nu+1}(t),b_n(t))$.
Recalling the definition of (not yet arranged) phase angles $\vartheta_i^1,\vartheta_j^2$,
we observe that
$$\vartheta_i^1\equiv \vartheta_{a_{i}}^1 \quad i\in\{1,\ldots,n-\nu\}
\hspace{4mm}\mbox{ and } \hspace{4mm}
\vartheta_j^2\equiv \vartheta_{-b_{j}}^2 \quad j\in \{n-\nu+1,\ldots, n\}.$$
\end{rem}

\noindent Let us now denote (cf. \cite{CodLev55}) by $\eta_j(a)$ the monotone sequence of the simple eigenvalues of the problem
\begin{equation}
\label{eig}
\left\{ \begin{array}{l}
u^{\prime \prime}(t)+ (a(t)\,+\,\eta)\,u(t)= 0, \\
\noalign{\medskip}
u(0)=0=u(1).
\end{array}\right.
\end{equation}
Recall that $\lim_{j \to +\infty} \eta_j(a) =+\infty$. Moreover, the eigenfunction corresponding to $\eta_j(a)$ has
exactly $(j-1)$ zeros on $(0,1)$. From Sturm's theory, a relation can be established between the 
eigenvalues $\eta_j(a)$ and the phase angle $\vartheta_{a+\eta}^1(t)$ associated to
\begin{equation}
\label{linsca}
\left\{\begin{array}{l}
u^{\prime}=v,\hspace{6mm} v'=-(a(t)\,+\,\eta)\,u,\\
\noalign{\medskip}
u(0)=0, \hspace{4mm} v(0)=1.
\end{array}\right.
\end{equation}
More precisely,
%\marginpar{NEW!!!}
\begin{eqnarray}
\label{eq:comp0}
\vartheta_{a+\eta}^1(1)=j\pi \hspace{6mm} \Longleftrightarrow \hspace{6mm}
\eta=\eta_j(a),\\
\label{comp1}
\vartheta_{a+\eta}^1(1)>j\pi \hspace{6mm} \Longleftrightarrow \hspace{6mm}
\eta>\eta_j(a),\\
\label{comp2}
\vartheta_{a+\eta}^1(1)<j\pi \hspace{6mm} \Longleftrightarrow \hspace{6mm}
\eta<\eta_j(a).
\end{eqnarray}
We remark that when $a$ is constant, it is possible to write
the explicit expression of $\eta_j(a)$ for each $j\in\N$. More precisely, if $a(t)= a\in\R$
for every $t\in[0,1]$, then
\begin{equation}
\label{eq:acost}
\eta_j(a):=j^2\pi^2-a\qquad \forall\,j\in\N.
\end{equation}

\bigskip
\bigskip
 
\section{Asymptotically linear Hamiltonian
systems}\label{sec:asympt}

\noindent We shall be concerned with the differential system
\begin{equation}\label{eq:secondorderequation}
J u''(t) + \nabla V(t,u(t)) =0,
\end{equation}
where $\nabla V(t,0)=0$ for all $t$. In what follows, we 
shall assume 
that $V: [0,1]\times \R^n \to \R$ is such that uniqueness and global continuability of solutions to
the initial value problems  associated to \eqref{eq:secondorderequation} are guaranteed. Since $u=0$ is a trivial solution of
\eqref{eq:secondorderequation}, it follows that the linearized
system at zero takes the form $ J u''(t) + D^2 V(t,0)u =0. $ By
setting
\begin{equation}\label{eq:Smatrice}
S(t, u) := \int_0^1 D^2V(t, su) ds, \qquad \forall\, (t,u) \in [0,1]
\times \R^n,
\end{equation}
system \eqref{eq:secondorderequation} can be written as follows
\begin{equation}\label{eq:secondosemplificato}
J u''(t) + S(t,u(t))u(t) =0.
\end{equation}
\noindent In this Section we describe the set of assumptions used for our main results; we then discuss some useful facts, in terms of Maslov index and phase angles, which hold in this framework. 

\noindent Let us first give the following
\begin{defin}\label{eq:sistemanondeg}
Given a path of symmetric matrices $S$, we say that the equation
$J u''(t) + S(t)u(t) =0$ is non-degenerate if the (linear) Dirichlet boundary
value problem 
\begin{equation}
\left\{\begin{array}{ll}
J u''(t) + S(t)u(t) =0, \quad  t\in [0,1],\\
u(0)=0=u(1)
\end{array}\right.
\end{equation}
has only the trivial solution.
\end{defin}
\noindent From now on we always assume that the strongly indefinite system
\eqref{eq:secondorderequation} is {\em
asymptotically linear at zero and infinity\/} meaning that there
exist two continuous paths of symmetric and split matrices $S_0$ and $S_\infty$ such that the following conditions hold:
\begin{enumerate}
\item[$(V_0)$] $S(t, \xi) \xi  = S_0(t) \xi + o(|\xi|)$ for $|\xi|\to 0$, uniformly in
$t$;
\item[$(V_\infty)$] $S(t, \xi)\xi= S_\infty(t) \xi + o(|\xi|)$ for $|\xi|\to +
\infty$, uniformly in $t$.
\end{enumerate}

\noindent Moreover, in what follows we suppose

\begin{enumerate}
\item[$(V_1)$] $S(t,x)$ is split for every $(t,x)\in [0,1] \times
\R^n$;
\item[$(V_2)$] the equations $J u''(t) + S_\infty(t)u(t) =0$
and $J u''(t) + S_0(t)u(t) =0$ are non-degenerate.
\end{enumerate}

\noindent The absence of $(V_2)$ would not affect the possibility of obtaining a multiplicity result like our Theorem \ref{eq:nondipdaalpha}; the only difference would then be in the exact count of the solutions with prescribed "nodal properties". The choice we make of using $(V_2)$ depends on the availability of useful formulas for the computation of the Maslov index (cf. \eqref{eq:formulaMaslovsplitdirichlet}, \eqref{eq:formulaMaslov}).

\noindent Finally, denoting by $W_i$ the $i$-th $(n-1)$-dimensional coordinate hyperplane in $\R^n$, we assume (as in \cite{CapDamPap}) that the following condition is fulfilled
\begin{enumerate}
\item[$(V_3)$] for every $i \in \n{n}$, the restriction of the
matrix $S$ to $[0,1] \times W_i$ is diagonal; i.e.
\[
S(t, x) := \diag(\lambda_1^i(t,x), \dots, \lambda_n^i(t,x)), \qquad
\forall\, (t, x)\in [0,1] \times W_i.
\]
\end{enumerate}

\noindent From $(V_3)$, it follows that $S_0$ and $S_\infty$ are diagonal matrices too.

\noindent Observe now that under the assumption $(V_\infty)$ there exists a
constant $M>0$ such that
\begin{equation}\label{eq:elastic property}
\ \ |\big(u(0),
u'(0)\big)|\leq R\Rightarrow |\big(u(t), u'(t)\big)|\leq
Re^{\max\{1,M\}}\ \textrm{for each}\ \ t \in [0,1].
\end{equation}
We refer to \cite{CapDamPap} for more details on this so-called "elastic property".

\noindent For the (nonlinear) self-adjoint second order boundary value
problem
\begin{equation}\label{eq:sistemasecondoordine}
\left\{\begin{array}{ll} J u''(t) + S(t,u(t))u(t) =0 \\
u(0)=0=u(1),
\end{array}\right.
\end{equation}
let us now consider for every $\alpha \in \R^n$ the Cauchy problem:
\begin{equation}\label{eq:problemadicauchylocale}
\left\{
\begin{array}{ll}
Ju''(t) + S(t, u(t))u(t)=0\\
u(0)=0,\quad J u'(0)=\alpha .
\end{array}
\right.
\end{equation}
Under the regularity assumptions on $S$,  we know
that, for each $\alpha\in \R^n$, there exists a unique solution $u_\alpha$ to
\eqref{eq:problemadicauchylocale}. Observe that $u_\alpha$ is a solution of
\eqref{eq:sistemasecondoordine} if and only if $u_\alpha(1) =0.$
\begin{defin}\label{def:Lsistem}
We call $\mathscr L$-system associated to
\eqref{eq:secondosemplificato} at $u_\alpha$ the linear second order
system
\begin{equation}\label{eq:Lsistema}
Ju''(t)+ S_\alpha(t) u(t)=0
\end{equation}
where $S_\alpha(t):= S(t, u_\alpha(t))$ for each $t \in [0,1]$.
\end{defin}
\noindent For every $\alpha \in \R^n$, we can develop the phase angle
analysis for the linear system \eqref{eq:Lsistema} and in particular
we can define, according to the previous notation, the matrices
$X_{\alpha,j}(t), X'_{\alpha,j}(t)$, the unitary symmetric matrices
$Y^1_\alpha :=(X'_{\alpha,1}(t)+ i X_{\alpha,
1}(t))(X'_{\alpha,1}(t)- i X_{\alpha, 1}(t))^{-1}$, $Y^2_\alpha
:=(X'_{\alpha,2}(t)+i X_{\alpha, 2}(t))(X'_{\alpha,2}(t)-i
X_{\alpha, 2}(t))^{-1}$, the angles $\vartheta_{l,\alpha}^1,\vartheta_{l,\alpha}^2$ with $l\in\n{n}$ and
the poset $(\Theta_\alpha,\prec^\oplus)$. Analogous definitions can be given for the "asymptotic linear systems"' at zero and infinity $Ju''(t)+ S_0(t) u(t)=0,$ and $Ju''(t)+ S_{\infty}(t) u(t)=0$, respectively.

\begin{defin}\label{def:angolidifaseconalpha}
We say that
$u_\alpha$ is a nontrivial $\underline h$-type solution if there
exists $\underline h \in (\N^{*n}, \poset)$ such that
$\theta^1_{l,\alpha}(1)= h_l\pi$ for all $l\in \{1,\ldots,n-\nu\}$ and
$\theta^2_{l,\alpha}(1)= -h_{l}\pi$ for all $l\in \{n-\nu+1,\ldots,n\}$,
where $\theta^j_{l,\alpha}$ are the phase angles associated to
\eqref{eq:Lsistema} for $j\in\n{2}$. In particular,  $(\N^{*n}, \poset)$ means
$(\N^{*(n-\nu)}, \preceq) \oplus (\N^{*\nu}, \preceq^{op})$.
\end{defin}

\noindent We can now state the following

\begin{lem}\label{thm:Hinfinito}
Under the assumptions $(V_0)-(V_\infty)$ we have
\begin{enumerate}
\item $S_\alpha$ tends to $S_\infty$ for $|\alpha| \to + \infty$
strongly in the $L^1$-norm topology;
\item $S_\alpha$ tends to $S_0$ for $|\alpha| \to  0$
strongly in the $L^1$-norm topology.
\end{enumerate}
\end{lem}
\proof We refer to \cite[Proposition 4.4,
Proposition 4.5]{CapDamPap}. \finedim 

\noindent As a direct consequence of
Lemma \ref{thm:Hinfinito} the result below easily follows.

\begin{cor}\label{thm:convergenzaflussieangoli}
Under the assumptions $(V_1)-(V_0)-(V_\infty)$ we have:
\begin{enumerate}
\item $Y^j_\alpha$ tends to $Y^j_\infty$ in the $C^0$-norm topology, for $|\alpha|
\to + \infty$;
\item $Y^j_\alpha$ tends to $Y^j_0$ in the $C^0$-norm topology for $|\alpha| \to 0$.
\end{enumerate}
\end{cor}

\noindent In what follows we shall associate an index to a linear second
order system of the form \eqref{eq:linearesecondoordine}; then, we will show how it can be
computed in some particular cases. 

\noindent Consider the fundamental solution $\phi$ of the first order system \eqref{eq:ham}. Being $\phi(\cdot)$ a path of symmetric matrices, it has a Maslov index $\mu_{L_0}\,(\phi,[\varepsilon ,1])$, where $\varepsilon$ is chosen in such a way that there are no crossings in $(0,\varepsilon]$; for brevity in what follows we shall write $m(S):=\mu_{L_0}\,(\phi,[\varepsilon ,1])$.

\begin{rem}\label{eps}
The existence of the constant $\varepsilon$ is guaranteed by the fact that the crossing instants of $\phi$
cannot accumulate at $0$ (cf. \cite{MPP}). Thus, the Maslov index $m(S)$ is
well defined and it is independent of the choice of $\varepsilon.$  
\end{rem}

\noindent In what follows we shall be concerned with some results on the computation of the Maslov index. 

\vspace{3mm}

\pparagraph{\noindent The Maslov index for constant and split matrices.} Consider
the second order Dirichlet boundary value problem
\begin{equation}\label{eq:seconODEconstant}
\left\{\begin{array}{lll} J u''(t)+ S u(t)=0 & &  \forall \, t \in [0,1]\\
u(0)=0=u(1).
\end{array}\right.
\end{equation}
\noindent For any real number $a$, let us consider the integer
\begin{equation}\label{eq:interomaslov}
N(a):= \#\{i \in \N^*| i^2 \pi^2 < a\}.
\end{equation}
Assume that \eqref{eq:seconODEconstant} has only the trivial solution. It is shown in \cite{MusPejPor07} that the following formula holds:
\begin{equation}\label{eq:formulaMaslovsplitdirichlet}
m(S) \, = \,  \sum_{i=1}^{n-\nu} N(\lambda_i) -
\sum_{i=1}^{\nu} N(-\mu_i),
\end{equation}
where $\lambda_i$ and $ \mu_i$ are the eigenvalues of $A$ and $B$,
respectively. In the particular case when $\nu =0$ the previous
formula reduces to:
\begin{equation}\label{eq:formulaMaslov}
m(S)  =  \sum_{i=1}^n N(\lambda_i).
\end{equation}
On the other hand, if the equation in \eqref{eq:seconODEconstant} is degenerate then, by Definition \ref{def:RobSalindex}, the following estimate holds:
\begin{equation}
\label{eq:stimeMaslovdegenere}
\left|\,m(S)-\sum_{i=1}^{n-\nu} N(\lambda_i) +
\sum_{i=1}^{\nu} N(-\mu_i)\,\right|\leq \frac{n}{2}.
\end{equation}

\begin{rem}\label{rmk:confrontocadapa03}
If $\nu=0$ and in the non-degenerate
case, the integer $m(S)$ agrees with
the total number of moments of verticality of  $S$ used
in \cite[Definition 3.6]{CapDamPap}. In particular, formula \eqref{eq:formulaMaslov} agrees with the formula given in \cite[Remark 3.9]{CapDamPap}. Observe that in the
degenerate case the integer used in \cite{CapDamPap}  does not coincide with the Maslov index; indeed if we compute (in the scalar case) the Maslov index from Definition \ref{def:RobSalindex}, an extra term $\pm 1/2$ appears.
\end{rem}

\vspace{3mm}

\noindent \pparagraph{The Maslov index for non constant and split
matrices.} In what follows, for brevity (and according to 
Remark \ref{eps}) we shall write the Maslov index $\mu_{\widetilde{L_0}}(\psi, [\varepsilon,1])$ of the fundamental solution $\psi$ of the
Hamiltonian system \eqref{eq:hamdisaccoppiato} (with respect to the
symplectic form $\widetilde \sigma:= \sigma_{n-\nu} \oplus
\sigma_{\nu}$) with 
$\mu_{\widetilde{L_0}}(\psi)$.

\begin{lem}\label{thm:3.13}
If $\psi(1)(\widetilde L_0)\cap \widetilde L_0 = \{0\}$, then the
Maslov index of $\psi$ is given by
\begin{equation}\label{eq:maslovconcondizionigenerali}
\mu_{\widetilde{L_0}}(\psi) \, = \,  \sum_{l=1}^{n-\nu}
k^1_l(1)-\sum_{l=1}^{\nu} k^2_l(1),
\end{equation}
where the integers $k_l^j$ have been defined in \eqref{eq:angolieinteri}.
\end{lem}
\proof Recall at first that
\[
\mu_{\widetilde{L_0}}( \psi) \, = \,  \mu_{L_0^A\oplus L_0^B }(
\psi).
\]
Denote by $\hat\psi_A$ the fundamental solution of the first order
system in $\R^{2(n-\nu)}$
\begin{equation}\label{primopezzo}
\left\{\begin{array}{ll} u'_1 = v_1  \\
v'_1= -A u_1
\end{array}\right.
\end{equation}
and by $\psi_B$ the fundamental solution of the first order system
in $\R^{2\nu}$
\begin{equation}\label{secondopezzo}
\left\{\begin{array}{ll} u'_2 = -v_2  \\
v'_2= -B u_2.
\end{array}\right.
\end{equation}
The direct sum property of the Maslov index implies that
\[
\mu_{\widetilde L_0}( \psi)  \,  =  \, \mu_{L_0^A }( \hat\psi_A) +
\mu_{L_0^B }( \psi_B).
\]
If we denote by $\hat \psi_{-B}$ the fundamental solution of the
first order system in $\R^{2\nu}$
\begin{equation}\label{terzopezzo}
\left\{\begin{array}{ll} u'_2 = z_2  \\
z'_2= B u_2,
\end{array}\right.
\end{equation}
then it follows that
\[
\mu_{L_0^A }( \hat\psi_A) + \mu_{L_0^B }( \psi_B) \, =  \, \mu_{L_0^A }(
\hat\psi_A) - \mu_{L_0^B }(\hat \psi_{-B}).
\]
Systems \eqref{primopezzo} and \eqref{terzopezzo} are equivalent to
the second order systems
$$
u''_1+Au_1=0, \;\;\;  u''_2-Bu_2=0.
$$
These systems are of the form studied in \cite{CapDamPap}. Now the
thesis follows by using \cite[Proposition 3.12]{CapDamPap}.\finedim
\begin{lem}\label{thm:strafava} The following equality holds
\begin{equation}\label{dueindici}
\mu_{\widetilde{L_0}}( \psi) \,=  \, m(S).
\end{equation}
\end{lem}
\proof Let us introduce on $\R^{2n}$ the symplectic
form $\omega_1$, by setting $\omega_1(z_1, z_2) := \langle
\widetilde\sigma z_1, z_2 \rangle$ for all $z_i \in \R^{2n}$ with
$i=1,2$. The proof then immediately follows by the
naturality property of the Maslov index (cf. \cite{RobSal95}),
combined with the fact that the matrix $U:(\R^{2n}, \omega)\to
(\R^{2n},\omega_1)$ is a symplectic isomorphism. \finedim

\noindent We end this section with some preliminary consequences of our assumptions.

\vspace{3mm}

\pparagraph{\noindent Maslov index and phase angles.} We denote by $m_0$ and
$m_\infty$, respectively, the Maslov indices of the fundamental
solution of the linear Hamiltonian systems at zero and at $\infty$.
(As a direct consequence of Lemma \ref{thm:strafava}, we do not need
to specify the symplectic structure we are referring to). Assuming
$m_0+n <m_\infty$, we define the following set
\begin{equation}
\label{eq:S}
\mathscr S:=\left\{\underline h \in (\N^{*n}, \poset) \colon \langle
\underline h, \underline j\rangle \in \big(m_0 + n-\nu,
m_\infty-\nu\big)\right\}.\
\end{equation}

\begin{rem}
\label{thm:sprimo}
In the case $m_\infty+n < m_0$, it is enough to define the set
$\mathscr S$ as follows
\begin{equation} \label{eq:S'} \mathscr S':=\left\{\underline h \in
(\N^{*n}, \poset) \colon \langle \underline h, \underline j\rangle
\in \big( m_\infty+n-\nu, m_0-\nu\big)\right\}.
\end{equation}
\end{rem}

\begin{lem}\label{thm:stimetraccecampinondef}
If $\mathscr S\not= \emptyset$ then there exists $\varepsilon >0 $
such that the following inequalities hold:
\begin{equation}\label{eq:disugtraccecampi}
\varepsilon <\langle\underline h, \pi\underline
j\rangle-(m_0+n-\nu)\pi , \qquad \textrm{and}\quad (m_\infty-\nu)
\pi-\langle\underline h, \pi\underline j\rangle>\varepsilon, \qquad
\forall\, \underline h \in \mathscr S.
\end{equation}
Moreover, there exists $\alpha_0:=\alpha_\varepsilon$
small enough such that
\begin{equation}
\label{eq:zzero} |\alpha|\leq \alpha_0\quad  \Rightarrow \qquad
\langle \Theta_\alpha(1), \underline 1 \rangle < \langle
\Theta_0(1), \underline 1\rangle + \varepsilon\quad
\textrm{and}\quad \langle \Theta_\alpha(1), \underline 1 \rangle <
(m_0+ n-\nu)\pi+ \varepsilon.
\end{equation}
Furthermore there exists $\alpha_\infty>\alpha_\varepsilon >0$ such
that
\begin{equation}
\label{eq:iinfty} |\alpha|\geq \alpha_\infty\quad  \Rightarrow
\qquad \langle \Theta_\alpha(1), \underline 1 \rangle >
\langle\Theta_\infty(1), \underline 1\rangle - \varepsilon\quad
\textrm{and}\quad \langle \Theta_\alpha(1), \underline 1 \rangle >
(m_\infty -\nu)\pi- \varepsilon.
\end{equation}
\end{lem}
\proof The conclusion follows from Corollary
\ref{thm:convergenzaflussieangoli} and  Lemma
\ref{thm:3.13}.\finedim 

\noindent Define $R:=\alpha_\infty$. By
\eqref{eq:elastic property} it follows that there exists $M>0$ such that
\begin{equation}\label{eq:proprietaelastica}
|\alpha| \leq \alpha_\infty \quad \Rightarrow \quad |(u_\alpha(t),
u'_\alpha(t))|\leq \alpha_ \infty e^{\max\{1,M\}} \qquad \forall\, t
\in [0,1].
\end{equation}
For $r \in (0,R)$, let $\mathscr D_r^R$ be the conical shell defined by
\[
\mathscr D_r^R := \{ x \in \R^{n} \colon \ \ r \leq |x| \leq R,\ \
x_i \geq 0, \ \ \forall\, i \in \n{n}\}.
\]
Recall that, for each $i \in \n{n}$, we denote by $W_i$ the $i$-th
$(n-1)$-dimensional coordinate hyperplane in the Euclidean space
$\R^n$. Let $\underline \alpha_i:= (\alpha_1,\dots, \alpha_{i-1}, 0,
\alpha_{i+1}, \dots, \alpha_n) \in \mathscr D_{\alpha_{0}}^{\alpha_{\infty}} \cap W_i$ and let
us consider the corresponding $\mathscr L$-system given by:
\[
Ju''(t)+ S_{\underline \alpha_i}(t) u(t)=0,
\]
where $S_{\underline\alpha_i}(t):= S(t, u_{\underline \alpha_i}(t))$
for the solution $u_{\underline \alpha_i}$ of the initial value
problem \eqref{eq:problemadicauchylocale} with $\alpha=\underline\alpha_i$. Consider also the eigenvalues $\lambda^{\underline \alpha_i}_1(t) , \dots ,
\lambda^{\underline \alpha_i}_n(t)$
of  $S_{\underline \alpha_i}(t)$.

\noindent Let
us denote by $D_i^{n-1}$ the $(n-1)$-dimensional closed disk of
radius $\alpha_\infty e^{\max\{1,M\}}$ contained in the hyperplane
$W_i$ and by $C_i$ the $n$-dimensional (full) cylinder $[0,1]
\times D_i^{n-1}$.
Since $\underline\alpha_i\in \mathscr D_{\alpha_{0}}^{\alpha_{\infty}} \cap W_i$, by \eqref{eq:proprietaelastica} it follows that
$u_{\underline \alpha_i}(t)\in D_i^{n-1}$ for each $t\in[0,1]$.
Thus, by assumption $(V_3)$
\begin{equation}
\label{eq:Salphai}
S_{\underline\alpha_i}(t) = \diag\left(\lambda_1^i(t, u_{\underline \alpha_i}(t)), \dots, \lambda_n^i(t, u_{\underline \alpha_i}(t))\right),
\qquad\forall\, t\in [0,1],
\end{equation}
whence we deduce that $\lambda^{\underline\alpha_i}_k(t)\equiv\lambda_k^i(t,u_{\underline \alpha_i}(t))$. Define
\begin{equation}\label{eq:autovalorimaggiorantipositivodefinito}
\overline \lambda_k^i:=\max\{\lambda_k^{i}(z): \,z\in C_i\}\qquad \forall\, {k=1, \dots, n}.
\end{equation}
Notice that
\begin{equation}
\label{eq:boh}
\lambda^{\underline\alpha_i}_k(t)\leq \overline\lambda_k^i \qquad \forall\,t\in[0,1].
\end{equation}
For each $i$, we define the two sets of  permutations
\begin{equation}
\label{eq:permumu}
\sigma^1_i: \{1, \dots, n-\nu\} \to \{1, \dots, n-\nu\}\qquad
\textrm{and}\quad  \sigma^2_i: \{n-\nu +1, \dots, n\} \to \{n-\nu
+1, \dots, n\}
\end{equation}
and the vectors
\begin{equation}\label{eq:alpostodellamatricecheserve}
\overline \Lambda_{\sigma^1_i}:=\big(\overline
\lambda_{\sigma^1_i(1)}^{i}, \dots,\overline
\lambda_{\sigma^1_i(n-\nu)}^{i}\big),\quad
\textrm{and}  \quad \overline
\Lambda_{\sigma^2_i}:=\big(\overline
\lambda_{\sigma^2_i(n-\nu +1)}^{i},
\dots,\overline \lambda_{\sigma^2_i(n)}^{i}\big)
\end{equation}
obtained by respectively arranging in increasing order the
components of
\[
\overline \Lambda_1^i:=\big(\overline
\lambda_1^{i}, \dots,\overline
\lambda_{n-\nu}^{i} \big),\ \,\,\, \textrm{and of} \ \,\,\,
\overline \Lambda_2^i:=\big(\overline
\lambda_{n-\nu+1}^{i}, \dots,\overline
\lambda_{n}^{i} \big).
\]
Denoting by $\overline \Delta$ the  diagonal matrix given by
\[
\diag\big(\overline \lambda^{1}_{\sigma^1_1(1)},
\dots ,\overline \lambda^{n-\nu}_{\sigma^1_{n-\nu}(n-\nu)},\overline
\lambda^{n-\nu+1}_{\sigma^2_{n-\nu+1}(n-\nu+1)},
\dots, \overline \lambda^{n}_{\sigma^2_n(n)}\big),
\]
we consider the second order boundary value problem
\begin{equation}\label{eq:sistemaausiliarioNONdirichletNondefinitopositivo}
\left\{\begin{array}{ll} Ju''(t) + \overline \Delta \, u(t) =0,
\qquad
t \in [0,1]\\
u(0)=0=u(1).
\end{array}
\right.
\end{equation}

\noindent For each $\underline h =(h_1, \dots, h_n) \in (\N^{*n}, \poset)$ we set (recalling that $\eta_{h_{i}}$ are the eigenvalues of problem \eqref{eig})
\begin{equation}
\label{eq:etah}
\overline\delta_{\underline{h}}\!:=\!\left(\eta_{h_{1}}\left(\overline\lambda^{1}_{\sigma^1_1(1)}\right)\!,
\dots,\eta_{h_{n-\nu}}
\left(\overline \lambda^{n-\nu}_{\sigma^1_{n-\nu}(n-\nu)}\right)\!,
\eta_{h_{n-\nu+1}}
\left(-\overline\lambda^{n-\nu+1}_{\sigma^2_{n-\nu+1}(n-\nu+1)}\right)\!,
\dots,
\eta_{h_{n}}
\left(-\overline \lambda^{n}_{\sigma^2_n(n)}\right)\!\right)\!,
\end{equation}
and we introduce the following set
\begin{equation}\label{eq:insiemeTNOnpositividefinito}
\mathscr T:=\left\{\underline h \in (\N^{*n}, \poset) \vert \quad
\overline\delta_{\underline h}\succeq^\ominus 0 ,\quad
m_0+n-\nu< \langle\underline h, \underline j\rangle<m_\infty -\nu
\right\}.
\end{equation}
In the next section we shall give sufficient conditions which guarantee that $\mathscr T$ is not empty. 
 
\noindent Note that, according to \eqref{eq:acost},
in the positive definite case $\nu=0$, $\mathscr T$
coincides with the set $\mathcal{T}$ defined in \cite{CapDamPap}.

\bigskip
\bigskip

\section{The main results}\label{sec:risultati}

\noindent The main idea in order to prove our results is to use the {\em
Miranda's fixed point theorem.\/} For the sake of completeness, we
recall it in a formulation suitable for the situation we are
dealing with.

\begin{thm}\label{thm:mirandatype}( \cite[Theorem 2.1]{CapDamPap} ).
Let $f: \mathscr D_r^R \to \R^{n} $ be a continuous vector field and
assume that the following conditions hold:
\begin{enumerate}
\item $\sum_{i=1}^{n} f_i(\alpha) <0 \ \ \textrm{for}\ \
|\alpha|=r$;\qquad  $\sum_{i=1}^{n} f_i(\alpha) >0 \ \ \textrm{for}\
\ |\alpha|=R$;
\item $f_i(\alpha) <0 \ \ \textrm{for}\ \ \alpha \in \mathscr D_r^R\cap
W_i$ \ \ \textrm{and}\,\,  $i \in \n{n}$.
\end{enumerate}
Then there exists at least one point $\widetilde \alpha$ in the
interior of $\mathscr D_r^R$ such that $f(\widetilde \alpha)=0$.
\end{thm}

\begin{rem}
\label{rem:disopplat}
The statement of Theorem \ref{thm:mirandatype} holds true if we replace
condition $(2)$ with
\begin{itemize}
\item[$(2')$]
$f_i(\alpha) >0 \ \ \textrm{for}\ \ \alpha \in \mathscr D_r^R\cap
W_i$ \ \ \textrm{and}\,\,  $i \in \n{n}$.
\end{itemize}
\end{rem}
\noindent Now, let $\mathscr D$ be the conical
shell $\mathscr D_{\alpha_0}^{\alpha_\infty}$; for any
$\underline h \in \mathscr S$, define $f\colon \mathscr D\to
\R^n$ as the continuous vector field whose components
are given by
\begin{equation}\label{eq:definizionedelcampo}
f_i(\alpha):=\left\{\begin{array}{ll}
\theta^1_{i, \alpha}(1) - h_i \pi , & \, i \in
\{1, \dots, n-\nu\}, \qquad \\
\noalign{\medskip} h_{i} \pi +\theta^2_{i, \alpha}(1) & \, i \in
\{n-\nu+1,\ldots,n\}
\end{array}\right.
\end{equation}
where, for $j \in {\bf 2}$, $\theta^j_{i, \alpha}$ are the phase angles
associated to the $\mathscr L$-system \eqref{eq:Lsistema}.

\begin{lem}\label{thm:primaMirandavalenondirichletnondefinito}
Assume $(V_1)-(V_2)-(V_0)-(V_\infty)$. Then the following inequalities hold:
\[
\sum_{i=1}^{n} f_i(\alpha) <0 \ \ \textrm{for}\ \
|\alpha|=\alpha_0;\qquad \sum_{i=1}^{n} f_i(\alpha) >0 \ \
\textrm{for}\ \ |\alpha|=\alpha_\infty.
\]
\end{lem}
\proof These are consequences of the second and
third inequalities in Lemma \ref{thm:stimetraccecampinondef} and
\eqref{eq:disugtraccecampi}. \finedim

\noindent We are now ready to state and prove our main result.

\begin{mainthm}
\label{eq:nondipdaalpha}
Let $n \geq 2$. Assume that the conditions
$(V_0)-(V_\infty)-(V_1)-(V_2)-(V_3)$ hold. Suppose that
\begin{equation}
\label{eq:nonvuoto}
\mathscr T\neq\emptyset.
\end{equation}
Then the boundary value problem \eqref{eq:sistemasecondoordine} has $2^n$ distinct $\underline h$-type solutions,
for every $\underline h\in\mathscr T$.
\end{mainthm} 

\proof
We fix $\underline h \in\mathscr T$
and we prove at first the
existence of $\widetilde\alpha=(\widetilde\alpha_1,\ldots,\widetilde\alpha_n)\in \mathscr D$
and of a solution $u$ of $\underline h$-type such that
$Ju'(0)=\widetilde\alpha$.
To this end, let $f\colon\mathscr D\to \R^n$ be the continuous vector field whose components
are defined in the equation \eqref{eq:definizionedelcampo}. By
taking into account Lemma
\ref{thm:primaMirandavalenondirichletnondefinito} it follows that
the first condition of Theorem \ref{thm:mirandatype} holds.

\noindent In order to conclude the proof of the Theorem it is enough to show
that also the second condition of  Theorem \ref{thm:mirandatype}
holds, i.e.
\begin{equation}
\label{eq:tesi}
\left\{\begin{array}{lll} f_i(\underline\alpha_i):=\theta^1_{i, \underline\alpha_{i}}(1) -
h_i \pi <0, & \quad \textrm{for each}\ \  i \in \{1,\ldots,n-\nu\},
& \underline\alpha_i \in \mathscr D \cap W_i
\\
f_i(\underline\alpha_i):=h_{i}\pi + \theta^2_{i, \underline\alpha_{i}}(1)<0, &
\quad \textrm{for each}\ \ i \in\{n-\nu+1,\ldots,n\}, & \underline\alpha_i \in \mathscr D \cap W_i.
\end{array}\right.
\end{equation}
Let us fix $\underline\alpha_i \in \mathscr D \cap W_i$.
As observed in \eqref{eq:Salphai}, by assumption $(V_3)$ it follows that for each $i\in\n{n}$
$$S_{\underline\alpha_i}(t)\diag\left(\lambda^{\underline\alpha_i}_1(t), \ldots, \lambda^{\underline\alpha_i}_n(t)\right).$$
Hence, taking into account Remark \ref{remrem} and the definition of phase angles for \eqref{eq:Lsistema}, we first note that
\begin{equation}
\label{eq:varte}
\vartheta^1_{k, \underline\alpha_{i}} \equiv \vartheta^1_{\lambda^{\underline \alpha_i}_{k}}
\hspace{5mm} \mbox{and} \hspace{5mm}
\vartheta^2_{k, \underline\alpha_{i}} \equiv \vartheta^2_{-\lambda^{\underline \alpha_i}_k}.
\end{equation}

\noindent First, we fix $i\in\{1,\ldots,n-\nu\}$. Observe that (due to the Sturm comparison principle stated in Lemma \ref{thm:favata}) the permutation $\sigma^1_i$
introduced in \eqref{eq:permumu} to arrange in increasing order the
constants $\overline\lambda_k^{i}$ arranges in increasing order the angles $\vartheta_{ \overline\lambda^{i}_{k} }^1(1)$ as well, 
i.e.
\begin{equation}
\label{eq:siconcludera1}
\vartheta_{ \overline\lambda^{i}_{\sigma^{1}_{i}(k)} }^1(1)\leq
\vartheta_{ \overline\lambda^{i}_{\sigma^{1}_{i}(h)} }^1(1)\qquad \mbox{ if }\,\,
k\leq h.
\end{equation}
On the other hand, in general the permutation $\sigma^1_i$ does not arrange the angles $\vartheta^1_{k, \underline\alpha_{i}}(1)$.
Indeed, in general, $\theta^1_{k, \underline\alpha_{i}}(1)\neq \vartheta^1_{\sigma^1_i(k), \underline\alpha_{i}}(1)$.
However, recalling that the definition of $\theta^1_{i, \underline\alpha_{i}}(1)$ comes from the arrangement
in increasing order of the angles $\vartheta^1_{k, \underline\alpha_{i}}(1)$, or, equivalently, of the angles
$\vartheta^1_{\sigma_i^1(k), \underline\alpha_{i}}(1)$, we infer that
\begin{equation}
\label{eq:terminera1}
\forall\,i\in\{1,\ldots,n-\nu\}\qquad \exists\,k_i\in\{1,\ldots,i\}:\qquad
\theta^1_{i, \underline\alpha_{i}}(1)\leq\vartheta^1_{\sigma^{1}_{i}(k_i), \underline\alpha_{i}}(1).
\end{equation}

\noindent Now, we fix $i\in\{n-\nu+1,\ldots,n\}$.
Taking into account that the permutation $\sigma^2_i$ introduced in \eqref{eq:permumu} arranges in increasing order the
constants $\overline\lambda_k^{i}$, by applying Lemma \ref{thm:favata} it is easy to verify that
\begin{equation}
\label{eq:siconcludera2}
\vartheta_{- \overline\lambda^{i}_{\sigma^{2}_{i}(k)} }^2(1)\leq
\vartheta_{- \overline\lambda^{i}_{\sigma^{2}_{i}(h)} }^2(1)\qquad \mbox{ if }\,\,
k\leq h.
\end{equation}
Moreover, recalling that also $\theta^2_{i, \underline\alpha_{i}}(1)$ comes from the arrangement
in increasing order of the angles $\vartheta^2_{\sigma_i^2(k), \underline\alpha_{i}}(1)$, we conclude that
\begin{equation}
\label{eq:terminera2}
\forall\,i\in\{n-\nu+1,\ldots,n\}\qquad \exists\,l_i\in\{n-\nu+1,\ldots,i\}:\qquad
\theta^2_{i, \underline\alpha_{i}}(1)\leq\vartheta^2_{\sigma^{2}_{i}(l_i), \underline\alpha_{i}}(1).
\end{equation}

\noindent As a next step, taking into account the relations \eqref{eq:boh} and \eqref{eq:varte},
we apply the Sturm comparison principle stated in Lemma \ref{thm:favata} to prove that
\begin{eqnarray}
\label{eq:acabara1}
\vartheta^1_{k, \underline\alpha_{i}}(1)
\leq\vartheta_{ \overline\lambda^{i}_{k} }^1(1),
&\quad & k\in\{1,\ldots,n-\nu\},\\
\label{eq:acabara2}
\vartheta^2_{k, \underline\alpha_{i}}(1)\leq \vartheta_{-\overline\lambda^{i}_{k} }^2(1),
&\quad & k\in\{n-\nu+1,\ldots,n\}.
\end{eqnarray}
Thus, by combining \eqref{eq:terminera1}, \eqref{eq:acabara1} and \eqref{eq:siconcludera1} with
the fact that $k_i\leq i$, we deduce that
\begin{equation}
\label{eq:finira1}
\theta^1_{i, \underline\alpha_{i}}(1)\leq\vartheta^1_{\sigma^{1}_{i}(k_i), \underline\alpha_{i}}(1)\leq
\vartheta_{ \overline\lambda^{i}_{\sigma^{1}_{i}(k_i)} }^1(1)\leq
\vartheta_{ \overline\lambda^{i}_{\sigma^{1}_{i}(i)} }^1(1),
\quad  i\in\{1,\ldots,n-\nu\}.
\end{equation}
To complete the proof, we recall that
$\eta_{h_i}\left(\overline\lambda^{i}_{\sigma^1_i(i)}\right)>0$
since $\underline h \in\mathscr T$. Hence, by \eqref{comp2}, we obtain
\[\vartheta_{\overline\lambda^{i}_{\sigma^{1}_{i}(i)}}^1(1)
<h_i\pi, \hspace{5mm} i \in \{1,\ldots,n-\nu\},\]
which implies the validity of the first inequalities in \eqref{eq:tesi}, i.e.
\begin{equation}
\label{eq:dig1}
f_i(\underline\alpha_i):=\theta^1_{i, \underline\alpha_{i}}(1) -  h_i \pi
<0,  \qquad \forall\,  i \in \{1, \dots, n- \nu\}.
\end{equation}
To prove the validity of the second inequalities of \eqref{eq:tesi},
we first combine \eqref{eq:terminera2}, \eqref{eq:acabara2} and \eqref{eq:siconcludera2} with
the fact that $l_i\leq i$ to infer
\[\theta^2_{i, \underline\alpha_{i}}(1)\leq\vartheta^2_{\sigma^{2}_{i}(l_i), \underline\alpha_{i}}(1)\leq
\vartheta_{ -\overline\lambda^{i}_{\sigma^{2}_{i}(l_i)} }^2(1)\leq
\vartheta_{ -\overline\lambda^{i}_{\sigma^{2}_{i}(i)} }^2(1),
\quad  i\in\{n-\nu+1,\ldots,n\}.\]
Secondly, recall that $\eta_{h_i}\left(-\overline\lambda^{i}_{\sigma^2_i(i)}\right)<0$
since $\underline h \in\mathscr T$.
Thus, from \eqref{comp1}, we deduce that
$$\vartheta_{-\overline\lambda^{i}_{\sigma^{2}_{i}(i)}}^1(1)
>h_i\pi, \hspace{5mm} i \in \{n-\nu+1,\ldots,n\},$$
which, according to Lemma \ref{thm:favata}, can be equivalently written as
\[\vartheta_{-\overline\lambda^{i}_{\sigma^{2}_{i}(i)}}^2(1)
<-h_i\pi, \hspace{5mm} i \in \{n-\nu+1,\ldots,n\}.\]
We can finally conclude that
$$f_i(\underline\alpha_i):= \theta^2_{i, \underline\alpha_{i}}(1)+h_{i} \pi
<0,  \quad \forall \,i \in \{n-\nu+1,\ldots,n\},$$
which completes the proof of \eqref{eq:tesi}.
Thus, Theorem \ref{thm:mirandatype} guarantees the existence
of $\widetilde \alpha$ in the interior of $\mathscr D$ such that
$f_i(\widetilde \alpha)=0$ for every $i\in\n{n}$.
Taking into account Lemma \ref{eq:angolidifasemolteplicita},
we deduce that {\em all the solutions of}
\begin{equation}
\label{eq:tildealpha}
J u''+S_{\widetilde\alpha} (t) u(t)=0,\qquad u(0)=0
\end{equation}
verify $u(1)=0$. Since $u_{\widetilde{\alpha}}$ solves
\eqref{eq:tildealpha}, we have proved the existence of a solution
$u$ of \eqref{eq:sistemasecondoordine} of $\underline h$-type such that
$J u'=\widetilde\alpha$ and $u'(0)\succeq^\ominus 0$.

\noindent In order to prove the existence of the other solutions it is enough
to apply the abstract result in the remaining $2^{n-1}$ conical
shells contained in the remaining hyper-octants determined by the
coordinate planes.\finedim
\begin{rem}
\label{rem:asimma}
Consider $\underline h=(h_1,\ldots,h_n) \in\mathscr T$.
By definition,
\begin{equation}
\label{eq:negativisigh!}
\eta_{h_i}\left(-\overline\lambda^{i}_{\sigma^2_i(i)}\right)<0,\qquad
\forall\, i\in\{n-\nu+1,\ldots,n\}.
\end{equation}
From \eqref{eq:acost},
we know that $\eta_{h_i}\left(-\overline\lambda^{i}_{\sigma^2_i(i)}\right)h_i^2\pi^2+\overline\lambda^{i}_{\sigma^2_i(i)}$,
and, consequently,
\[\overline\lambda^{i}_{\sigma^2_i(i)}<-h_i^2\pi^2\leq -\pi^2\qquad \forall\,i\in\{n-\nu+1,\ldots,n\}.\]
In particular, we have shown that $\overline\lambda^{i}_{\sigma^2_i(i)}$ should be negative
to guarantee that $\mathscr T\neq \emptyset$.
\end{rem}
\begin{rem}
We may redefine $\overline\delta_{\underline h}$ and ${\mathscr T}$ by replacing in
\eqref{eq:etah} and \eqref{eq:insiemeTNOnpositividefinito}
the constants $\overline \lambda_k^i$ introduced in \eqref{eq:autovalorimaggiorantipositivodefinito}
with the following functions
\begin{equation}
\label{eq:tempo}
\overline \lambda_k^i(t):=\max\{\lambda_k^{i}(t,x): \,x\in D_i\}\qquad \forall\, k,i\in\n{n}.
\end{equation}
Then, it is easy to prove that (under minor modifications) Theorem \ref{eq:nondipdaalpha} holds true in this more general setting.\\
We remark that the maps $\overline\lambda^{i}_{\sigma^2_i(i)}(\cdot)$,
obtained by arranging in increasing order the maps \eqref{eq:tempo},
do not need to be negative in the
whole interval $[0,1]$ to satisfy \eqref{eq:negativisigh!}; it is sufficient that they are negative on a subset
of $[0,1]$ of positive measure.
\end{rem}

\noindent Now we give a sufficient condition in order to guarantee
that $\mathscr T\neq \emptyset$.
\begin{lem}\label{lem:suffcond}
A sufficient condition in order that $\mathscr T\neq \emptyset$
is given by
\begin{equation}
\label{eq:condsuffnondef}
\begin{array}{l}
\displaystyle{\overline\lambda^{i}_{\sigma^2_i(i)}<-\pi^2\hspace{9mm} \forall\,i\in\{n-\nu+1,\ldots,n\}},\\
\noalign{\bigskip}
\displaystyle{m_0-\frac{n}{2}+\nu<m(\overline \Delta)<m_\infty-\frac{5}{2} n+\nu.}
\end{array}
\end{equation}
\end{lem}

\proof
Denoting by $m(\overline \Delta)$ the Maslov index associated to the system
\eqref{eq:sistemaausiliarioNONdirichletNondefinitopositivo}, by
formula \eqref{eq:stimeMaslovdegenere}
it directly follows that
\begin{equation}
\label{eq:45}
m(\overline \Delta)-\frac{n}{2}\leq \sum_{i =1}^{n-\nu} N\left(\overline
\lambda^{i}_{\sigma_i^1(i)}\right) -\sum_{i=n-\nu+1}^n
N\left(-\overline \lambda^{i}_{\sigma_i^2(i)}\right)\leq m(\overline \Delta)+\frac{n}{2}.
\end{equation}
According to \eqref{eq:acost}, we easily observe that
\begin{eqnarray}
\label{eq:boh1}
N\left(\overline\lambda^{i}_{\sigma^1_i(i)}\right)<h_i-1
\hspace{3mm} &\Longrightarrow& \hspace{3mm}
\eta_{h_i}\left(\overline\lambda^{i}_{\sigma^1_i(i)}\right)>0,
\hspace{9mm} i \in \{1,\ldots,n-\nu\},\\
\label{eq:boh2}
N\left(-\overline \lambda^{i}_{\sigma_i^2(i)}\right)\geq h_i
\hspace{5mm} &\Longrightarrow& \hspace{3mm}
\eta_{h_i}\left(-\overline\lambda^{i}_{\sigma^2_i(i)}\right)<0,
\hspace{5mm} i \in \{n-\nu+1,\ldots,n\},
\end{eqnarray}
for each $h_l\in\N$, with $l\in\n{n}$.
Let us define
\begin{eqnarray*}
&\tilde h_i:=N\left(\overline\lambda^{i}_{\sigma^1_i(i)}\right)+2,  \qquad &
i \in \{1,\ldots,n-\nu\}\\
&\tilde h_i:=N\left(-\overline \lambda^{i}_{\sigma_i^2(i)}\right),  \qquad &
i \in \{n-\nu+1,\ldots,n\}.
\end{eqnarray*}
According to the first assumption in \eqref{eq:condsuffnondef},
$\tilde h_i\geq 1$ for each $i \in\n{n}$.

\noindent Our aim consists in proving that the vector
$\underline{\widetilde{h}}:=(\tilde h_1,\ldots,\tilde h_n)$
belongs to $\mathscr T$.\\
From \eqref{eq:boh1}-\eqref{eq:boh2}, we immediately deduce that
\[\overline\delta_{\underline{\widetilde{h}}}\succeq^\ominus 0.\]
Moreover, by formula \eqref{eq:45}, we obtain that
\[m(\overline \Delta)-\frac{n}{2}+2(n-\nu)\leq
\langle\underline{\widetilde{h}}, \underline j\rangle\leq
m(\overline \Delta)+\frac{n}{2}+2(n-\nu).\]
Hence, taking into account the second assumption in \eqref{eq:condsuffnondef}, we infer
that
\[m_0+n-\nu=m_0-\frac{n}{2}+\nu-\frac{n}{2}+2(n-\nu)<
\langle\underline{\widetilde{h}}, \underline j\rangle
<m_\infty-\frac{5}{2} n+\nu+\frac{n}{2} +2(n-\nu)=m_\infty-\nu.\]
Recalling the definition of $\mathscr T$, we conclude that
$\underline{\widetilde{h}}\in\mathscr T$.
\finedim

\noindent Note that the first condition in \eqref{eq:condsuffnondef}
agrees with Remark \ref{rem:asimma}.

\noindent According to Lemma \ref{lem:suffcond}, we can thus write
\begin{mainthm}
\label{thm:CaDaPo}
Let $n \geq 2$. We assume that the conditions
$(V_0)-(V_\infty)-(V_1)-(V_2)-(V_3)$ hold.
Suppose moreover that $m_0+2n<m_\infty$.
If condition \eqref{eq:condsuffnondef} is satisfied,
then the boundary value problem \eqref{eq:sistemasecondoordine} has $2^n$ distinct $\underline h$-type solutions,
for every $\underline h\in\mathscr T$.
\end{mainthm}

\noindent According to formula \eqref{eq:formulaMaslovsplitdirichlet}, we note
that condition \eqref{eq:condsuffnondef} can be refined if
the equation in \eqref{eq:sistemaausiliarioNONdirichletNondefinitopositivo}
is non-degenerate.
More precisely, we can prove the following result.
\begin{mainthm}
\label{thm:3}
Let $n \geq 2$. We assume that the conditions
$(V_0)-(V_\infty)-(V_1)-(V_2)-(V_3)$ hold.
Moreover, suppose that the equation $J u''(t)+\overline \Delta u(t)=0$
is non-degenerate, and that the following conditions are satisfied:
\begin{equation}
\label{eq:nondegi}
\begin{array}{l}
\displaystyle{\overline\lambda^{i}_{\sigma^2_i(i)}<-\pi^2\hspace{9mm} \forall\,i\in\{n-\nu+1,\ldots,n\}},\\
\noalign{\bigskip}
\displaystyle{m_0<m(\overline \Delta)<m_\infty- n.}
\end{array}
\end{equation}
Then the boundary value problem \eqref{eq:sistemasecondoordine} has $2^n$ distinct $\underline h$-type solutions,
for every $\underline h\in\mathscr T$.
\end{mainthm}
\proof
To prove this result, we can argue exactly as before.
We have to use formula \eqref{eq:formulaMaslovsplitdirichlet} instead of
\eqref{eq:stimeMaslovdegenere} and we should
take into account that in the non-degenerate case the relation \eqref{eq:boh1}
can be relaxed into the following
\begin{equation}
\label{eq:ciao}
N\left(\overline\lambda^{i}_{\sigma^1_i(i)}\right)<h_i
\hspace{3mm} \Longrightarrow \hspace{3mm}
\eta_{h_i}\left(\overline\lambda^{i}_{\sigma^1_i(i)}\right)>0,
\hspace{9mm} i \in \{1,\ldots,n-\nu\},
\end{equation}
since, for each $i\in\{1,\ldots,n-\nu\}$, there is no $k\in\N$
such that $\overline\lambda^{i}_{\sigma^1_i(i)}=k^2\pi^2$.
According to \eqref{eq:ciao}, we should now choose
$\tilde h_i:=N\left(\overline\lambda^{i}_{\sigma^1_i(i)}\right)+1$, for
each $i \in \{1,\ldots,n-\nu\}$. Then, the thesis easily follows.
\finedim
 
\begin{rem}
A result analogue to Theorem \ref{eq:nondipdaalpha} can be written by reversing the first inequality in the definition of the set 
$\mathscr T$; moreover, according to Remark \ref{thm:sprimo}, an analogous statement can be written by exchanging $m_0$ and $m_{\infty}$ in the second inequality in the definition of the set 
$\mathscr T$. In both cases, variants of Theorem \ref{thm:CaDaPo} and Theorem \ref{thm:3} can be obtained accordingly.
 \end{rem}

\begin{rem}\label{notcomp}
Note that \eqref{eq:sistemasecondoordine} can be written in the
equivalent form
\begin{equation}
\label{eq:CDP}
\left\{\begin{array}{ll} u''(t) +J S(t,u(t))u(t) =0 \\
u(0)=0=u(1).
\end{array}\right.
\end{equation}
The split assumption $(V_1)$ guarantees that
$JS:[0,1]\to\Bsym(\R^n)$ remains a continuous path of symmetric
matrices, whenever $S:[0,1]\to\Bsym(\R^n)$ is  continuous. In
particular, the presence of this symmetry enables us to handle
problem \eqref{eq:CDP} with the methods employed in
\cite{CapDamPap}. Given the equation
\[u''(t) +J S(t)u(t) =0,\] and introduced $z:=(u,u')$,
one can denote by $\hat{\psi}_{(JS)}$ the fundamental solution of
\[z'= \begin{pmatrix} 0 & {\rm Id} \cr -JS & 0
\end{pmatrix}z.\]
Assuming the validity of conditions
$(V_0)-(V_\infty)-(V_1)-(V_2)-(V_3)$ and applying direclty Theorem
4.7 in \cite{CapDamPap}, it is possible to prove the existence of
$2^n$ nontrivial $\underline h$-type solutions, whenever $\underline h$
belongs to a suitable, non empty subset $\hat{\mathscr T}$ of
\[
\hat{\mathscr S}:=\left\{\underline h \in (\N^{*n}, \prec)
\colon \langle \underline h, \underline 1\rangle \in \big(\mu_0 + n,
\mu_\infty\big)\right\}\
\]
where
$\mu_0$, $\mu_\infty$ are, respectively, the
Maslov indices of the fundamental solutions $\hat{\psi}_{(JS_0)}$
and $\hat{\psi}_{(JS_\infty)}$ relative to $L_0$.

\noindent According to Lemma \ref{thm:3.13} and to the previous
notation, we can observe that
\[\begin{array}{lll}
m_0 := \mu_{L_0^{A_0} }( \hat\psi_{A_0}) -\mu_{L_0^{B_0} }(
\hat\psi_{-B_0})& while &
\mu_0:=\mu_{L_0^{A_0} }(\hat\psi_{A_0}) +\mu_{L_0^{B_0} }( \hat\psi_{-B_0})\\
m_\infty := \mu_{L_\infty^{A_\infty} }( \hat\psi_{A_\infty}) -
\mu_{L_\infty^{B_\infty} }(\hat\psi_{-B_\infty})& while &
\mu_\infty:=\mu_{L_\infty^{A_\infty} }(\hat\psi_{A_\infty}) +
\mu_{L_\infty^{B_\infty} }( \hat\psi_{-B_\infty})
\end{array}
\]
whenever
\[
S_0 =\begin{pmatrix} A_0 & 0 \\
0 & B_0
\end{pmatrix}
\quad and \quad
S_\infty =\begin{pmatrix} A_\infty & 0 \\
0 & B_\infty
\end{pmatrix}.
\]
These relations allow us to conclude that the two approaches are
deeply different, since, in general, they provide solutions with different nodal properties. Indeed, if we focus our attention on the sets $\mathscr S$
$($whose definition is given in \eqref{eq:S}$)$ and $\hat{\mathscr
S}$, we notice that the intervals $\big(m_0 +
n-\nu,m_\infty-\nu\big)$ and $\big(\mu_0 + n,\mu_\infty\big)$
coincide in the case $\nu=0$, but, in general, they are not comparable.
\end{rem}

\noindent We end this section with two propositions which deal with the question of the emptiness/nonemptiness of $\mathscr T$.

\noindent We first consider the case when we are in presence of radial symmetry.
\begin{prop}
\label{asimmetria}
\begin{equation}
S(t,x)=S(t,|x|)\; {\rm for \, every} \; (t,x)\in[0,1] \times \R^n \qquad \Longrightarrow \qquad \mathscr T=\emptyset.
\end{equation}
\end{prop}
\proof
According to assumption $(V_3)$ and definitions
\eqref{eq:autovalorimaggiorantipositivodefinito}-\eqref{eq:permumu}, we observe that
$\overline \lambda_k^i$ and, consequently, $\sigma^1_i,\,\sigma^2_i$ are independent of $i\in\n{n}$.
In particular, let us set
%\marginpar{NEW!}
\begin{equation}
\label{eq:modulo}
\overline \lambda_k^i:=\overline\lambda_k,\qquad \sigma^j_i:=\sigma^j,\qquad
\forall i,\,k\in\n{n},\hspace{2mm} \forall j\in\n{2}.
\end{equation}
Arguing by contradiction, assume the existence of $\underline h\in\mathscr T$.\\
Fix $\underline\alpha_1\in\mathscr D \cap W_1$.
As in \eqref{eq:terminera1}, we observe that for every $k\in\{1,\ldots,n-\nu\}$
there exists $l_k\in\{1,\ldots,k\}$ such that
$\theta^1_{k, \underline\alpha_{1}}(1)\leq\vartheta^1_{\sigma^{1}(l_k), \underline\alpha_{1}}(1)$.
With the same argument used to achieve \eqref{eq:finira1},
by combining \eqref{eq:acabara1} and \eqref{eq:siconcludera1} with
the fact that $l_k\leq k$, we deduce that
%\begin{equation}
%\label{eq:finira3}
\[\theta^1_{k, \underline\alpha_{1}}(1)\leq\vartheta^1_{\sigma^{1}(l_k), \underline\alpha_{1}}(1)\leq
\vartheta_{ \overline\lambda_{\sigma^{1}(l_k)} }^1(1)\leq
\vartheta_{ \overline\lambda_{\sigma^{1}(k)} }^1(1),
\quad  \forall\,k\in\{1,\ldots,n-\nu\}.\]
%\end{equation}
Since $\underline h\in\mathscr T$, by definition,
\[\eta_{h_k}\left(\overline\lambda_{\sigma^1(k)}\right)>0,
\quad  \forall\,k\in\{1,\ldots,n-\nu\},\]
which, according to \eqref{comp2}, leads to
\[\vartheta_{\overline\lambda_{\sigma^{1}(k)}}^1(1)
<h_k\pi, \hspace{5mm} k \in \{1,\ldots,n-\nu\}.\]
Taking into account the definition of $f$ given in
\eqref{eq:definizionedelcampo}, we can conclude that
\[f_k(\underline\alpha_1):=\theta^1_{k, \underline\alpha_{1}}(1) -  h_k \pi
<0,  \qquad \forall\,  k \in \{1, \dots, n- \nu\},\hspace{2mm}
\forall \underline\alpha_1\in\mathscr D\cap W_1.\]
Following analogous steps, we can prove that
\[f_k(\underline\alpha_1):=\theta^2_{k, \underline\alpha_{1}}(1) + h_k \pi
<0,  \qquad \forall\,  k \in \{n- \nu+1,\ldots,n\},\hspace{2mm}
\forall \underline\alpha_1\in\mathscr D\cap W_1.\]
Considering $\alpha^*:=(0,\beta)$ with $\beta\in\R^{n-1}$ and
$|\beta|=\alpha_\infty$,
we observe that $\alpha^* \in\mathscr D\cap W_1$.
From the previous inequalities, we infer that
\[\sum_{k=1}^{n} f_k(\alpha^*) <0,\qquad |\alpha^*|=\alpha_\infty,\]
which contradicts Lemma \ref{thm:primaMirandavalenondirichletnondefinito}.
This implies that $\mathscr T= \emptyset$.
\finedim

\begin{rem}
\label{rem:mdeltaminfty}
In the radial symmetric case when $S(t,x)=S(t,|x|)$ for every $(t,x)\in[0,1]\times\R^n$, it is possible to prove that 
\begin{equation}
\label{eq:simmental}
m(\overline \Delta)>m_\infty-\frac{3}{2} n-\frac{\varepsilon}{\pi}.
\end{equation}
Moreover, if the equation $J u''(t)+\overline \Delta u(t)=0$ is non-degenerate, then
\eqref{eq:simmental} can be refined into the following
\[m(\overline \Delta)>m_\infty- n-\frac{\varepsilon}{\pi}.\]
The above inequalities show that in the radial case the sufficient condition \eqref{eq:condsuffnondef} for the non-emptiness of the set $\mathscr T$ is violated. Notice that in the case $n=\nu$, the contradiction follows from the fact that
$m(\overline \Delta)$ is a semi-integer, and, without loss of generality,
$\varepsilon$ in Lemma \ref{thm:stimetraccecampinondef} can be chosen
smaller than $\frac{\pi}{2}$.
\end{rem}

\noindent In order to end our discussion on the emptiness/non emptiness 
of the set $\mathscr T$, let us first observe that from Corollary \ref{thm:convergenzaflussieangoli} there exists $\tilde \alpha_{\infty}$ such that 
\begin{equation}\label{eq:contradiction}
|\alpha|\geq \tilde\alpha_\infty
\quad
\Rightarrow
\qquad
\theta^j_{i, \alpha}(1)> \theta^j_{i, \infty}(1)-\frac{\varepsilon}{n}
\qquad \forall\, i\in\n{n},\,\,\,\forall j\in \n{2},
\end{equation}
with $\varepsilon$ as in Lemma \ref{thm:stimetraccecampinondef}.

\noindent Remark also that \eqref{eq:contradiction} implies 
\begin{equation}
\label{eq:iiinfty} |\alpha|\geq \tilde\alpha_\infty\quad  \Rightarrow
\qquad \langle \Theta_\alpha(1), \underline 1 \rangle >
\langle\Theta_\infty(1), \underline 1\rangle - \varepsilon\quad
\textrm{and}\quad \langle \Theta_\alpha(1), \underline 1 \rangle >
(m_\infty -\nu)\pi- \varepsilon.
\end{equation}

\noindent On the other hand, if we go back to the constant $\alpha_\infty$ (whose existence is proved in Lemma \ref{thm:stimetraccecampinondef}) then we have 

\begin{prop}
\label{sceltaalphainfty}
\begin{equation}
\label{vuoto}
\alpha_\infty \, \geq \tilde\alpha_\infty  \qquad  \Rightarrow \qquad
\mathscr T = \emptyset.
\end{equation}
\end{prop}
\proof
Assume, by contradiction, the existence of $\underline h\in\mathscr T$.
For each $i\in\n{n}$ and for every vector $\beta^i\in \mathscr D\cap W_i$,
let us define
\[\Theta^*(\beta^1,\ldots,\beta^n):=\left(
\theta^1_{1, \beta^{1}}(1),\ldots, \theta^1_{n-\nu, \beta^{n-\nu}}(1),
\theta^2_{n-\nu+1, \beta^{n-\nu+1}}(1),\ldots,\theta^2_{n, \beta^{n}}(1)
\right).\]
Note that we can choose $|\beta^i|=\alpha_\infty$ for each $i\in\n{n}.$
During the proof of Theorem \ref{eq:nondipdaalpha},
we have proved the validity of inequalities \eqref{eq:tesi}, which lead to
\[\langle \Theta^*(\beta^1,\ldots,\beta^n), \underline 1\rangle
<\langle\pi\underline h, \underline j\rangle,
\qquad \forall\,\beta^i\in \mathscr D\cap W_i,\hspace{3mm} i\in\n{n}.\]
Moreover, by combining \eqref{eq:contradiction}
with the fact that $\alpha_\infty \, \geq \tilde\alpha_\infty$, we get
\begin{equation}
\label{eq:versolameta}
\langle \Theta^*(\beta^1,\ldots,\beta^n), \underline 1\rangle>
\langle \Theta_\infty(1), \underline 1\rangle-\varepsilon>
(m_\infty-\nu)\pi-\varepsilon\qquad \mbox{if }
|\beta^i|=\alpha_\infty\quad \forall\,i\,\in\n{n},
\end{equation}
from which we conclude that
\[\langle\pi\underline h, \underline j\rangle>
(m_\infty-\nu)\pi-\varepsilon,\]
which contradicts \eqref{eq:disugtraccecampi}.
The emptiness of $\mathscr T$ follows.
\finedim

\noindent Thus, we have learnt that the non-emptiness of $\mathscr T$ is possible when  \eqref{eq:iinfty} (which deals with the whole vector $\Theta_{\alpha}$) is not a consequence of \eqref{eq:contradiction}. In other words, we are implicitely requiring that
\begin{equation}
\label{implic}
\alpha_\infty  \qquad   <  \qquad \tilde\alpha_\infty.
\end{equation}

\noindent According to Proposition \ref{asimmetria}, condition \eqref{implic} can be interpreted as the requirement that the radial symmetry of $Ju''+S_\infty(t)u(t)=0$ must not be preserved away from infinity.

\end{document}